\documentclass[a4paper,11pt]{amsart}

\usepackage{amsmath,graphics}
\usepackage{amssymb}
\usepackage{amsfonts}
\usepackage{latexsym}
\usepackage{eucal}
\usepackage{array}
\usepackage{enumerate}
\newtheorem{thm}{Theorem}[section]

\newtheorem{propo}[thm]{Proposition}
\newtheorem{lem}[thm]{Lemma}

\newtheorem{conj}[thm]{Conjecture}

\renewcommand{\Re}{{\rm Re}}
\renewcommand{\Im}{{\rm Im}}
\newcommand{\R}{\mathbb{R}}
\newcommand{\C}{\mathbb{C}}
\newcommand{\Z}{\mathbb{Z}}
\newcommand{\N}{\mathbb{N}}

\newcommand{\D}{ \Omega}
\newcommand{\B}{{\mathcal B}}

\newcommand{\G}{\mathbf{G}}
\newcommand{\GL}{\mathrm{GL}}
\newcommand{\End}{\mathrm{End}}

\newcommand{\lt}{{\mathcal L}}

\newcommand{\U}{{\mathcal U}}
\newcommand{\half}{{\textstyle{\frac{1}{2}}}}

\newcommand{\LG}{\mathfrak g }

\begin{document}
\bibliographystyle{plain}
\title[Heat Kernel and compact extensions]{Heat kernel
 and Rates of mixing for compact extensions of expanding maps}

\author{Fr\'ed\'eric Naud }
\address{
Laboratoire d'Analyse non-lin\'eaire et G\'eom\'etrie\\
Universit\'e d'Avignon, 33 rue Louis Pasteur\\
84000 Avignon France.
}
\email{frederic.naud@univ-avignon.fr}
\thanks{This paper is dedicated to the memory of Serge Lang who believed in the heat kernel as "a universal gadget which is a dominant factor
practically everywhere in mathematics".}

\begin{abstract}
 We consider skew-extensions $\widehat{T}:I\times  \G$ of real analytic markov expanding maps
 of the interval $T:I\rightarrow I$ by a compact Lie group $\G$. 
 We show that for the natural $\widehat{T}$-invariant measures associated to an analytic potential $\varphi$, 
 there exists an infinite dimensional cone of real valued and analytic observables for which the rate of 
 mixing enjoys an {\it explicit exponential lower bound} related to the topological pressure $P(2\varphi)$. 
 This estimate is consistent with a conjecture of Dolgopyat and Pollicott on the spectrum of complex transfer operators. 
\end{abstract}

 \maketitle

\begin{section}{Introduction}
Partially hyperbolic systems form a non structurally stable class of dynamical systems whose ergodic and statistical properties are
considered difficult to study, mainly because of the presence of a neutral bundle which prevents application of all classical tools designed for uniformly hyperbolic systems. Compact Lie groups extensions of uniformly hyperbolic systems is
an interesting familly of such maps where issues such as stable ergodicity and mixing can be studied in detail with the techniques of Lie groups and representation theory. For a detailed account on the history of the qualitative ergodic theory of such systems, we refer
the reader to the introduction of \cite{Dolgisrael} and references herein.
In this paper we will focus
on partially hyperbolic maps obtained as skew extensions of expanding maps of the interval. 
Let $T:I\rightarrow I$ be
\footnote{We will provide a more detailed definition in the next $\S$.}
 a Markov, topologically mixing, smooth enough, uniformly expanding map of the unit interval $I=[0,1]$. Let $\G$ be a compact connected Lie group endowed with the normalized Haar measure $m$, and consider $\tau:I\rightarrow \G$ a smooth map. The skew extension
 $\widehat{T}:I\times \G\rightarrow I\times \G$ of $T$ is defined by
 $$\widehat{T}(x,g)=(Tx,\tau(x)g).$$
 Let $\mu_\varphi$ be the unique equilibrium $T$-invariant measure related to a H\"older potential $\varphi$ on $I$,
 then the product measure $\widehat{\mu}_\varphi:=\mu_\varphi \times m$ is a $\widehat{T}$-invariant measure. A popular choice of $\varphi$
 would be $\varphi=-\log\vert T' \vert$ which yields the absolutly continuous SRB measure. Let $F,G$ be some
 reasonably smooth functions on $I\times \G$, then a fundamental question of ergodic theory is to investigate the asymptotic behaviour
 of the correlation function defined for $n\in \N$ by
 $$C_\varphi (F,G)(n):=\int_{I\times \G}(F\circ \widehat{T}^n) G d\widehat{\mu}_\varphi-
 \int_{I\times \G}Fd\widehat{\mu}_\varphi \int_{I\times \G}G d\widehat{\mu}_\varphi.$$
 Correlation functions are the basic objects that give a measure of the "chaoticity" of the system and are not only interesting
 for themselves but also because of their relationship with fine statistical properties such as the central limit theorem.
 One of the very first quantitative result is the exponential decay for {\it equivariant observables} that can be found in Field-Melbourne-T\"or\"ok \cite{Melbourne}. At the same time, Dolgopyat showed in the case of the SRB measure \cite{Dolgisrael} that for generic $\tau$ and H\"older observables $F,G$, 
 the correlation function decays exponentially fast as $n\rightarrow +\infty$. Moreover,
 this property is equivalent to stable ergodicity (within the set of skew extensions). The exponential bound obtained
 through Dolgopyat technique is however hardly explicit and does not tell us how the different parameters (lyapunov exponents, topological entropy) of the system
 influence this rate of mixing. On the other hand, still in the case of the SRB measure, it is very likely that the more recent machinery of Tsujii et al, \cite{AGT,Tsujii} should provide an explicit upper bound related to the maximal expansion
 rate of the map $\widehat{T}$. In \cite{Dolgisrael}, Dolgopyat raises the question of finding examples where {\it effective estimates} on correlations can be proved. In the case of non uniformly hyperbolic systems, {\it lower bounds} on the rate of mixing are notoriously hard to prove, and there exists to our knowledge no general results in that direction for exponentially mixing systems. Motivated by the lower bounds we obtained
 in \cite{Naud5} for real analytic semi-flows which are toy models of Anosov flows, we show in this paper that for {\it general
 real analytic compact extensions } $\widehat{T}$, it is possible to build an infinite dimensional manifold of real analytic observables for which correlations functions enjoy an {\it effective lower bound} related to the topological pressure of $2\varphi$.

 There is therefore, like in the case of semi-flows \cite{Naud5}, 
 a natural barrier to the speed of mixing which is related to entropy.
 We point out that this behaviour contrasts sharply with the case of uniformly hyperbolic systems where the rate of
 mixing for analytic observables can be superexponential, for example in the case of linear hyperbolic maps on the torus.
 
 \bigskip
 \noindent Let us give some more precise statements. First we will work with extensions by a compact connected Lie group $\G$, which is non-trivial i.e. $\G\neq \{Id\}$. The skew function $\tau:[0,1]\rightarrow \G$ is assumed to be real analytic \footnote{We will give a precise meaning to this in the next
 $\S$.}.
 About the expanding map $T:I\rightarrow I$, we will assume that 
$$I=[0,1]=\cup_{i=1}^k I_i ,\ k\geq 2, $$
where each $I_i=[a_i,b_i]$ is a closed interval and ${\rm Int}(I_i)\cap {\rm Int}(I_j)=\emptyset$ if $i\neq j$. 
For all $i=1,\ldots,k$, the map $T_i:=T\vert_{[a_i,b_i]}$ is {\it real analytic},  and is a one-to one and onto map $T_i:I_i\rightarrow [0,1]$ which is 
{\it uniformly expanding } i.e. 
$$\min_{x\in I_i}\vert T_i'(x) \vert >1.$$
Under these standard assumptions, $T:[0,1\rightarrow [0,1]$ is topologically mixing and is semi conjugated to the full shift 
$$\sigma: \{1,\ldots,k \}^\N \rightarrow  \{1,\ldots,k \}^\N $$
on $k$ symbols
\footnote{In addition, the intersection of two different intervals $I_i$ and $I_j$, if not empty, cannot contain periodic points of $T$ (only pre-periodic points). Therefore there is a one-to-one correspondence between periodic points of the shift map $\sigma$ and periodic points of $T$.}.
We will be concerned with $T$-Invariant equilibrium measures $\mu_\varphi$ where the potential $\varphi:I\rightarrow \R$ is also real analytic on 
$[0,1]$. The associated measure $\mu_\varphi$
is then the unique measure maximizing the topological pressure (variational formula) :
$$P(\varphi)=\sup_{\mu \in {\mathcal M}_{inv}} \left \{     h_\mu(T)+\int_0^1 \varphi d\mu \right \},$$
where ${\mathcal M}_{inv}$ is the set of $T$-invariant probability measures and $h_\mu(T)$ is the measure theoretic entropy. For simplicity, we will assume that $P(\varphi)=0$ which is the case for all popular potentials.

Under the above assumptions and notations, we prove the following result.
\begin{thm}
\label{mainthm}
There exists a constant $\gamma(\G)>0$
\footnote{The constant $\gamma(\G)$ depends essentially on the rank of $\G$ that is the dimension of maximal torii, an explicit
expression will be provided in $\S 7$.} such that for all $\epsilon>0$, one can find an infinite dimensional cone 
${\mathcal C}_\epsilon$ of real valued, analytic functions such that for all $F \in {\mathcal C}_\epsilon$,
$$\limsup_{n\rightarrow +\infty}\vert C_\varphi(F,F)(n)\vert^{\frac{1}{n}} \geq e^{\gamma(\G)P(2\varphi)-\epsilon}.$$
\end{thm}
\noindent {\bf Remarks.}
We always have the lower bound 
$$P(2\varphi)\geq \int_I \varphi d\mu_\varphi=-h_{\mu_\varphi}(T)\geq -h_{\rm top}(T)=-\log(k),$$
with equality if and only if $\varphi$ is cohomologous to $-\log(k)$ the potential of the measure of maximal entropy.
On the other hand in the case of the smooth invariant SRB measure
($\varphi=-\log\vert T' \vert$), it tells us that if the derivative of $T$ at a fixed point $x_0$ is of size $1+\varepsilon$, then $$P(2\varphi)\geq h_{\delta_{x_0}}(T)-2\int_I \log\vert T' \vert d\delta_{x_0}=-2\log(1+\varepsilon),$$ 
where $\delta_{x_0}$ is the dirac mass at $x_0$. Thus $P(-2\log \vert T' \vert)$ tends to $0$ as $\varepsilon \rightarrow 0$, hence exhibiting arbitrarily slow rates of (exponential mixing) for the extension
$\widehat{T}$. In the case of the measure of maximal entropy, our lower bound remains unaffected by the expanding properties of $T$. This is not surprising in view of the recent result of Go\"uezel \cite{Gouezel} who was able to obtain exponential rates of mixing (for the measure of maximal entropy) in the case of skew-extensions of intermittent markov maps with neutral fixed points.

\bigskip \noindent
Although the statements above bear some obvious similarities with the results we obtained for semi-flows, the proof requires some new ideas to deal with the compact factor $\G$ in full generality. In particular, we need to prove some a priori bounds on the growth of certain dynamical zeta functions that do not follow from
\cite{Naud5}. The key a priori estimates are derived from Sogge's bound \cite{Sogge} and Donnelly-Fefferman work \cite{Fefferman},
which are both sharp results on elliptic partial differential operators on compact manifolds. 
The core of the argument in \cite{Naud5} is a Selberg's like trace formula for semi-flows. In the present setting, we need to deal with a discrete family of transfer operator indexed by irreducible representations and we have to replace the trace formula argument by a more delicate analysis based on the heat kernel and volume estimates around maximal Tori.  
We believe, based on some heuristics on transfer operators \cite{DP} (see $\S 2$) that the term
$P(2\varphi)$ is {\it optimal}, whereas the constant $\gamma(\G)$ may not be.
We also point out that we have chosen in this paper to keep the hyperbolic factor as simple as possible, concentrating on the compact group factor $\G$. It is however clear from our technique that similar results can be obtained without major modifications for extensions of real analytic expanding map of the torus or real analytic Anosov diffeomorphisms by using the approach of Faure-Roy \cite{Faureroy} for example. 

\bigskip \noindent
The plan of the paper is as follows. In $\S$ 2 we prove the main result by reduction to a spectral result on a family of transfer operators twisted by the irreducible representations of the group $\G$.  In $\S 3$, we give the proof in the simplest case
of the one dimensional torus which is a good outline of the general case, without the difficulties of the abstract general setting.
Section $\S 4, \S 5$ are devoted to some a priori bounds on certain Fredholm determinants. These bounds are not obvious in general and require some elliptic partial differential equation arguments. In $\S 6$ we recall some basic facts on the heat kernel on $\G$, our main tool. Finally, we gather our results and prove the key theorem in $\S 7$ following the line of proof of the torus case. 
\end{section}
\begin{section}{Representations and twisted transfer operators}
 In this section, we show how to reduce Theorem \ref{mainthm} to a spectral result on a family of transfer operators. We use the notations of the introduction. Let $\varphi:I\rightarrow \R$ be an analytic potential with zero topological pressure $P(\varphi)=0$.
 The (positive) Ruelle transfer operator is defined as usual on continuous functions by 
 $$\lt_\varphi(f)(x)=\sum_{Ty=x} e^{\varphi(y)}f(y).$$ 
 The Classical Ruelle-Perron-Frobenius Theorem, see for example \cite{Baladi,PP}, states that there exists a unique probablity
 measure $\nu_\varphi$ such that for all $f \in C^0(I)$,
 $$\int_I f d\nu_\varphi=\int_I \lt_\varphi(f) d\nu_\varphi.$$
 In addition there exists an invariant positive continuous non-vanishing density $h_\varphi \in C^0(I)$ such that
 $$\lt_\varphi(h_\varphi)=h_\varphi.$$
 The equilibrium $T$-invariant probability measure $\mu_\varphi$ that realizes the pressure is, up to normalization, 
 $\mu_\varphi=h_\varphi \nu_\varphi$. 
 
 \bigskip \noindent
 Let $\pi:\G\rightarrow \GL(V)$ be a finite dimensional {\it irreducible} complex \footnote{In this paper, all the representations are assumed to be
 complex.} representation of $\G$, that is to say a continuous group
 homomorphism from $\G$ to the group of linear isomorphism of the complex vector space $V$ such that the left action
 of $\G$ on $V$ has no other invariant subspaces than $V$ itself and $\{ 0\}$. We assume in addition that $\pi$
 is non trivial, i.e. is not identically equal to the identity. Such a representation does exist since $\G$ is assumed to be non trivial.
 
 \bigskip \noindent 
 Let $A:I\rightarrow \End(V)$ be a continuous map from $I=[0,1]$ to the complex vector space of linear maps
 from $V$ to itself. For all $(x,g)\in I\times \G$, we set 
 $$F(x,g)=\mathrm{Tr}(A(x)\pi(g)),$$
 and let $H(x,g)\in C^0(I\times \G)$. Let us compute the correlation function $C_\varphi(H,F)(n)$.
 Because $\pi$ is irreducible and non trivial, we have by Fubini's theorem
 $$\int_{I\times \G} F d\widehat{\mu}_\varphi=\mathrm{Tr}\left \{ \int_I A(x) d\mu_\varphi \int_\G \pi(g) dm(g) \right \}=0.$$
 Therefore we have 
 $$C_\varphi(H,F)(n)=\int_{I\times \G} H(T^nx,\tau(T^{n-1}x)\ldots\tau(Tx) \tau(x) g) F(x,g) d\widehat{\mu}_\varphi(x,g).$$
 For any $f:I\rightarrow \G$, we 
 set $f^{[n]}(x):=f(x)f(Tx)\ldots f(T^{n-1}x)$. Using Left invariance of Haar measure we get
 $$C_\varphi(H,F)(n)=\int_{I\times \G} H(T^nx,g) F(x,(\tau^{-1})^{[n]}(x)g) d\widehat{\mu}_\varphi(x,g).$$
Using Fubini again and the action of the transfer operator $\lt_\varphi^n$, we obtain
$$C_\varphi(H,F)(n)=$$
$$\int_\G\int_{I} H(x,g)  h_\varphi(x)^{-1}\lt^n_\varphi \left \{ h_\varphi(y) F(y,(\tau^{-1})^{[n]}(y)g) \right \}(x) d\mu_\varphi(x)dm(g).$$
On the other hand we have 
$$\lt^n_\varphi \left \{ h_\varphi(y) F(y,(\tau^{-1})^{[n]}(y)g) \right \}(x)$$
$$=\mathrm{Tr}\left \{ 
\pi(g)h_\varphi^{-1}(x)\sum_{T^ny=x}e^{\varphi^{(n)}(y)} h_\varphi(y) A(y) \pi( (\tau^{-1})^{[n]}(y))\right\},$$
where $\varphi^{(n)}(x)=\varphi(x)+\ldots +\varphi(T^{n-1}x)$. It can be rewritten as
$$\lt^n_\varphi \left \{ h_\varphi(y) F(y,(\tau^{-1})^{[n]}(y)g) \right \}(x)
=\mathrm{Tr}\left \{ \pi(g) {\mathcal M}^n_\pi(A)(x)  \right\},$$
where ${\mathcal M}_\pi: C^0(I,\End(V))\rightarrow C^0(I,\End(V))$ is the linear operator defined by
$${\mathcal M}_\pi(A)(x):=h_\varphi^{-1}(x)\sum_{Ty=x}e^{\varphi(y)}h_\varphi(y) A(y) \pi(\tau^{-1}(y)).$$
To summarize the above computation, we have obtained
\begin{equation}
\label{cor}
C_\varphi(H,F)(n)=\int_{I\times \G}H(x,g) \mathrm{Tr}\{ \pi(g) {\mathcal M}_\pi^n(A)(x) \} d\widehat{\mu}_\varphi(x,g).
\end{equation}
Theorem \ref{mainthm} now follows from the next key result.
\begin{thm}
\label{key}
Under the assumptions of $\S1$, for all $\epsilon>0$, there exist infinitely many non equivalent irreducible finite dimensional complex 
representations $\pi_\alpha, \alpha \in {\mathcal S}_\epsilon$
$$\pi_\alpha:\G\rightarrow \GL(V_\alpha),$$ 
such that each operator ${\mathcal M}_{\pi_\alpha}$ has a real analytic eigenfunction $A_\alpha\not \equiv 0$
$${\mathcal M}_{\pi_\alpha}(A_\alpha)=\lambda_\alpha A_\alpha,$$
with the lower bound $\vert \lambda_\alpha \vert \geq e^{\gamma(\G)(P(2\varphi)-\epsilon)}$.
\end{thm}
This result shows the existence of non trivial point spectrum with a uniform lower bound for an infinite
family of transfer operators. In the simplest case of the one dimensional torus $\G=\R/2\pi\Z$,
irreducible representations $(\pi_q)_{q\in \Z}$ are one dimensional and given by the usual Fourier analysis i.e.
$\pi_q(\theta)=e^{-iq\theta}$ so that the transfer operator ${\mathcal M}_{q}$ is given by
$${\mathcal M}_{q}(A)(x):=h_\varphi^{-1}(x)\sum_{Ty=x}e^{\varphi(y)}h_\varphi(y) e^{iq\tau(y)}A(y).$$
This is the usual "complex" transfer operator which plays a key role in Dolgopyat's approach 
\cite{Dolgopyat} and subsequent works on hyperbolic flows \cite{BaladiVallee,PS,Polldo1,Naud2,Stoyanov}. 
A standard "diagonal approximation argument" in the spirit of quantum physics, see \cite{DP}, tells us that the following conjecture should hold.
\begin{conj}
For a generic choice of $\tau$, for all $\epsilon>0$, there exists $q_0 \geq 1$ such that for all
$\vert q \vert \geq q_0$, the spectral radius $\rho_{\textrm{sp}}$ of ${\mathcal M}_q$, when acting on smooth enough function spaces, satisfies
$$e^{\half P(2\varphi)-\epsilon}\leq \rho_{\textrm{sp}}({\mathcal M}_q) \leq e^{\half P(2\varphi)+\epsilon}.$$
\end{conj}
\noindent
In the one dimensional case $\gamma(\G)=3/2$, so up to a factor $3$, our lower bound agrees with that conjecture.

\bigskip
\noindent We now explain how to deduce the statement of Theorem \ref{mainthm} from the above result.
Let us recall a basic fact about conjugate representations. If $\pi:\G\rightarrow \GL(V_\pi)$ is an irreducible representation on a finite dimensional complex vector space $V_\pi$, there exists a Hermitian inner product on
$V_\pi$ denoted by $\langle .,.\rangle_\pi$ such that $\pi$ is unitary. It induces a natural antilinear isomorphism
$${\mathcal T}:\left \{  \begin{array}{ccc}
 V_\pi&\rightarrow &V'_\pi \\
  v&\mapsto& \langle .,v\rangle_\pi 
\end{array}\right.,$$
where $V'_\pi$ is the dual space.
For simplicity, if $v\in V_\pi$ we set ${\mathcal T}(v)=\overline{v}$. The dual space $V'_\pi$ is now endowed with
a hermitian inner product by setting 
$$\langle \overline{u},\overline{v} \rangle_{\overline{\pi}}:=\langle v,u\rangle_\pi.$$
The {\it conjugate representation} $\overline{\pi}:\G\rightarrow \GL(V')$ is defined by
$\overline{\pi}={\mathcal T}\circ \pi \circ {\mathcal T}^{-1}$. One can check that $\overline{\pi}$ is irreducible
and that for all $u,v \in V$ and $g\in \G$ we have
$$\overline{\langle \pi(g)u,v \rangle}_\pi=\langle \overline{\pi}(g)\overline{u},\overline{v}\rangle_{\overline{\pi}}.$$
Given $A \in \GL(V_\pi)$ we can define $\overline{A} \in \GL(V'_\pi)$ by $\overline{A} \overline{v}:=\overline{Av}$, and the following formula holds
\begin{equation}
\label{conjugate}
\overline{\textrm{Tr}(\pi(g) A)}=\textrm{ Tr}(\overline{\pi}(g) \overline{A}).
\end{equation}
In terms of characters (the character is the trace of the representation), we have the formula
$$\overline{\chi_\pi}(g)=\chi_{\overline{\pi}}(g)=\chi_\pi(g^{-1}). $$
Fix $\epsilon >0$ and let ${\mathcal S}_\epsilon$ be a (countable non finite) set of irreducible non equivalent representations as given by Theorem \ref{key}. By removing some elements of 
$$\{ \pi_\alpha, \alpha\in {\mathcal S}_\epsilon \},$$
 we can assume that for all $\alpha \in {\mathcal S}_\epsilon$, 
$\pi_\alpha$ is non trivial and either $\overline{\pi_\alpha}$ is equivalent to $\pi_\alpha$, $\overline{\pi_\alpha}$ is not equivalent to any $\pi_\beta$
for all $\beta \in {\mathcal S}_\epsilon$. Given $A_\alpha(x)$ a real analytic eigenfunction of 
${\mathcal M}_{\pi_\alpha}$ as given by Theorem \ref{key}, we set 
$$F_\alpha(x,g)=\Re(\textrm{Tr}\{\pi_\alpha(g) A_\alpha(x)\}),\ 
G_\alpha(x,g)=\Im(\textrm{Tr}\{\pi_\alpha(g) A_\alpha(x)\}).$$
Using Schur's classical orthogonality relations, we have
$$\int_{I\times \G} \vert \textrm{Tr}\{\pi_\alpha(g) A_\alpha(x)\} \vert^2 d\widehat{\mu}_\varphi=
\frac{1}{d_\pi}\int_I \textrm{Tr}\{ A_\alpha(x) A^*_\alpha(x) \} d\mu_\varphi(x),$$
where $d_\pi$ is the dimension of the representation and $A^*$ is the adjoint operator for the inner product
associated to $\pi$. Recall that the equilibrium measure $\mu_\varphi$ has a Gibbs property and is thus positive on open sets, so the above integral is non vanishing. We can therefore assume (by replacing $A_\alpha$ by 
$iA_\alpha$ if necessary) that for all $\alpha$, 
$$\int_{I\times \G} F_\alpha^2 d\widehat{\mu}_\varphi>0.$$
We will need the following Lemma.
\begin{lem}
Let $\pi_\alpha,\ \pi_\beta$ be two representations as above with $\alpha\neq \beta$, then we have
$$\int_{I\times \G} F_\alpha F_\beta d\widehat{\mu}_\varphi=
\int_{I\times \G} F_\alpha G_\beta d\widehat{\mu}_\varphi=0.$$
\end{lem}
\noindent {\it Proof}. Write 
$$F_\alpha=
\frac{1}{2}(\textrm{Tr}\{\pi_\alpha(g) A_\alpha(x)\}+ \overline{\textrm{Tr}\{\pi_\alpha(g) A_\alpha(x)\}}),$$
$$G_\alpha=
\frac{1}{2i}(\textrm{Tr}\{\pi_\alpha(g) A_\alpha(x)\}-\overline{\textrm{Tr}\{\pi_\alpha(g) A_\alpha(x)\}}),$$
and use Schur orthogonality of non-equivalent representations together with formula (\ref{conjugate}) on
conjugate representations. $\square$

\bigskip
\noindent
Let ${\mathcal A}_\epsilon$ be the real vector space generated by the family of functions $F\alpha$
$${\mathcal A}_\epsilon=\textrm{Span}\{ F_\alpha,\ \alpha \in {\mathcal S}_\alpha \}.$$
Pick $F\in {\mathcal C}_\epsilon:={\mathcal A}_\epsilon \setminus \{ 0\}$ and write 
$$F=\sum_{\alpha \in {\mathcal F}} a_\alpha F_\alpha,$$
where ${\mathcal F}$ is a non-empty finite subset of ${\mathcal S}_\alpha$, and $a_\alpha\neq 0$ for all $\alpha \in {\mathcal F}$. Computing the correlation function $C_\varphi(F,F)(n)$ using formula (\ref{cor}) yields
$$C_\varphi(F,F)(n)=\sum_{\alpha,\beta \in {\mathcal F}}a_\alpha a_\beta 
\int_{I\times \G} F_\alpha \Re\left \{\lambda_\beta^n \textrm{Tr}(\pi_\beta(g) A_\beta(x)) \right \}
d\widehat{\mu}_\varphi,$$
which using the orthogonality relations above turns out to be
\begin{equation}
\label{cor2}
C_\varphi(F,F)(n)=\sum_{\alpha \in {\mathcal F}}a_\alpha^2 \left( \rho_\alpha^n \cos(n \theta_\alpha)
\int F_\alpha^2 d\widehat{\mu}_\varphi -\rho_\alpha^n \sin(n\theta_\alpha)\int F_\alpha G_\alpha 
d\widehat{\mu}_\varphi \right),
\end{equation}
where we have set $\lambda_\alpha=\rho_\alpha e^{i\theta_\alpha}$, with $\rho_\alpha >0$ and $\theta_\alpha \in \R$. We can now conclude the proof by using the Dirichlet box principle in the following form.
\begin{lem}
 Let $\alpha_1,\ldots,\alpha_N \in \R$ and $D\in \N\setminus \{0\}$. For all $Q\geq 2$, one can find
 an integer $q \in \{D,\ldots ,DQ^N \}$ such that 
 $$\max_{1\leq j \leq N} \textrm{dist}(q\alpha_j, \Z) \leq \frac{1}{Q}.$$
\end{lem}
\noindent
The proof is an elementary use of the box principle and we omit it. We can apply
it to $\alpha_j=\theta_\alpha/2\pi$ to construct a sequence of integers $n_q$ such that
$n_q \rightarrow +\infty$ as $q\rightarrow +\infty$ with the property that for all $q$, for all 
$\alpha \in {\mathcal F}$, $$\cos(n_q\theta_\alpha)\geq \frac{3}{4}\ \textrm{and}\  \vert \sin(n_q \theta_\alpha)\vert\int \vert F_\alpha G_\alpha \vert d\widehat{\mu}_{\varphi}\leq \frac{1}{4} 
\int F_\alpha^2d\widehat{\mu}_{\varphi}.$$
Going back to formula (\ref{cor2}), we have the lower bound for all $q$
$$C_\varphi(F,F)(n_q)\geq  e^{n_q\gamma(\G)(P(2\varphi)-\epsilon)}\frac{1}{2} \int F^2d\widehat{\mu}_{\varphi},$$
which implies the statement of Theorem \ref{mainthm}
$$\limsup_{n\rightarrow +\infty}\vert C_\varphi(F,F)(n)\vert^{\frac{1}{n}} \geq e^{\gamma(\G)(P(2\varphi)-\epsilon)}.$$
\end{section}

\begin{section}{The simplest case of skew extensions by the torus}
Before we move on to the general case of extensions by a compact Lie group $\G$, we would like to explain the basic ideas
of the proof in the simplest case where $\G=\R/2\pi \Z$. This section will serve as a guideline for the more difficult case of an abstract compact connected Lie group. 
\begin{subsection}{Transfer operators on Bergmann spaces}
We need first to be more precise about the choice of function space on which the transfer operators ${\mathcal M}_q$ act.
For all $j=1,\ldots,k$, we set $\gamma_j(x)=T_j^{-1}(x)$ to be the inverse branches of $T$. Because $T$ is real analytic and
expanding, it is possible to show (see Bandtlow-Jenkinson \cite{BJ1}) that there exists $\varepsilon_0>0$ such that for
all $\varepsilon\leq\varepsilon_0$, the complex domain 
$$\Omega_\varepsilon:=\{z\in \C\ :\ \min_{x \in [0,1]} \vert z-x\vert <\varepsilon \}$$
is (strictly) contracted by each inverse branch $\gamma_j$. More precisely, for all $j=1,\ldots,k$, $\gamma_j$ extends
holomorphically to $\Omega_\varepsilon$ and has a continuous extension to $\overline{\Omega_\varepsilon}$, in addition
$$\gamma_j(\overline{\Omega_\varepsilon})\subset \Omega_\varepsilon.$$
By taking $\varepsilon$ small enough, we can assume that both $\varphi$ and $\tau$ have an holomorphic extension to 
$\Omega_\varepsilon$ and a continuous extension to $\overline{\Omega_\varepsilon}$.  
The classical Bergmann space $B^2(\Omega_\varepsilon)$ is defined as 
$$B^2(\Omega_\varepsilon):=\left \{ f:\Omega_\varepsilon\rightarrow \C,\:\ f\ \textrm{holomorphic}\ 
\textrm{and}\ \int_{\Omega_\varepsilon}\vert f \vert^2 dxdy<+\infty \right \}.$$
It is a Hilbert space when endowed with the norm 
$$\Vert f\Vert_{B^2}^2:=\int_{\Omega_\varepsilon}\vert f \vert^2 dxdy.$$ 
The transfer operators ${\mathcal M}_q:B^2(\Omega_\varepsilon)\rightarrow B^2(\Omega_\varepsilon)$ are defined by
$${\mathcal M}_q(f)(z)=h_\varphi^{-1}(z)\sum_{j=1}^k e^{\varphi(\gamma_j z)+iq\tau(\gamma_j z)}h_\varphi(\gamma_j z) f(\gamma_j z).$$
We remark that the above expression makes sense (provided that $\varepsilon$ is small enough) since it can be shown \cite{Naud5} that the density $h_\varphi$ which
is given by Ruelle-Perron-Frobenius Theorem is actually non-vanishing and analytic on $[0,1]$. Below we summarize the
spectral properties of ${\mathcal M}_q$ that we use in the sequel.
\begin{propo} Using the above notations we have:
\label{spec1}
\begin{enumerate}
\item Each ${\mathcal M}_q:B^2(\Omega_\varepsilon)\rightarrow B^2(\Omega_\varepsilon)$ is a compact trace class operator.
\item The spectral radius of ${\mathcal M}_q$ is bounded by $1$ with respect to $q$.
\item Let $\vert \lambda_0(q)\vert\geq \vert \lambda_1(q)\vert \geq \ldots \geq \vert \lambda_n(q)\vert$ denote the
sequence of eigenvalues of ${\mathcal M}_q$. There exists a constant $C>0$ and $0<\rho<1$ such that for all
$q\in \Z$ and all $n\in \N$,
$$\vert \lambda_n(q)\vert \leq Ce^{C\vert q \vert} \rho^n.$$
\item For all $q\in \Z$ and all $n\in \N$, we have the formula for the trace of ${\mathcal M}_q^n$
$$\textrm{Tr}({\mathcal M}_q^n)=\sum_{T^nx=x} e^{iq\tau^{(n)}(x)} \frac{e^{\varphi^{(n)}(x)}}{1-((T^n)'(x))^{-1}},$$
where the sum runs over $n$-periodic points of $T$ and $\tau^{(n)}(x)=\tau(x)+\tau(Tx)+\ldots+\tau(T^{n-1}x)$.
\end{enumerate}
\end{propo}
\noindent
Let us make some comments on the proof of this standard facts. Result $(1)$ dates back to the work of Ruelle on expanding maps
\cite{Ruelle0}. Bound $(2)$ follows from classical contraction bounds and the formula for the topological pressure, see \cite{Naud5}, Proposition 2.1. Bound $(3)$ follows \cite{BJ2}, Theorem 5.13. Notice that the critical factor $e^{C\vert q \vert}$
comes from estimating
$$\sup_{z\in \Omega_\varepsilon} \vert e^{iq\tau(z)} \vert.$$ 
These estimates, though they are trivial in the specific case of the torus, will require some special treatment in the full
generality of compact Lie groups. Formula $(4)$ is a standard computation of the trace that can be found for example in
\cite{BJ2}. 
\end{subsection}
\begin{subsection}{A priori bounds on determinants}
Just like in \cite{Naud5}, a priori bounds on dynamical determinants is a key step in the proof of Theorem \ref{key}. Since
we want to prove the existence of non-trivial spectrum for the family ${\mathcal M}_q$, it is natural to define for $\zeta \in \C$
$${\mathcal Z}_q(\zeta):=\det(I-\zeta {\mathcal M}_q),$$
which is an entire function of $\zeta$ and equals  the infinite product
$$ {\mathcal Z}_q(\zeta)=\prod_{n\in \N}(1-\zeta \lambda_n(q) ).$$
We will need the following result.
\begin{propo}
\label{bound1}
Assume that for all $\vert q\vert$ large, ${\mathcal Z}_q(\zeta)$ does not vanish on the open disc $D(0,\rho_1)$, for some $\rho_1>0$. Then for all $r<\rho_1$, one can find $C_r>0$ such that for all $\vert q \vert$ large and all $\zeta \in D(0,r)$,
$$\left \vert  \frac{{\mathcal Z}'_q(\zeta)}{{\mathcal Z}_q(\zeta)} \right \vert\leq  C_r \vert q \vert.$$
\end{propo}
\noindent {\it Proof}. For all $\zeta \in \C$, we can write
$$\vert {\mathcal Z}_q(\zeta) \vert \leq \exp \left(   \sum_{n=0}^{+\infty} 
\log(1+\vert \zeta \vert \vert \lambda_n(q)\vert )\right ). $$
By splitting the above sum and using estimates $(2)$ and $(3)$ from Proposition \ref{spec1}, we get 
$$\sum_{n=0}^{+\infty} 
\log(1+\vert \zeta \vert \vert \lambda_n(q)\vert )= \sum_{n=0}^N \log(1+\vert \zeta \vert \vert \lambda_n(q)\vert )
+\sum_{n=N+1}^{+\infty} \log(1+\vert \zeta \vert \vert \lambda_n(q)\vert )$$
$$\leq (N+1)\log(1+\vert \zeta \vert)+C\frac{\rho^{N+1}}{1-\rho}\vert \zeta \vert e^{C\vert q \vert}.$$
Choosing
$$N=\left [ \frac{C\vert q \vert}{\vert \log \rho \vert} \right ]+1,$$
we have for well chosen constants $C_1,C_2$ and large $\vert q \vert$
$$\log\vert {\mathcal Z}_q(\zeta) \vert \leq C_1\vert q \vert \log(1+\vert \zeta \vert)+C_2 \vert \zeta \vert.$$
Assume now that ${\mathcal Z}_q(\zeta)$ does not vanish on $D(0,\rho_1)$. Since ${\mathcal Z}_q(0)=1$, we can 
define a complex holomorphic logarithm, denoted by $\log{\mathcal Z}_q(\zeta)$ for all $\zeta \in D(0,\rho_1)$, with the property that $\log {\mathcal Z}_q(0)=0$ and 
$$\Re(\log {\mathcal Z}_q(\zeta) )=\log \vert {\mathcal Z}_q(\zeta) \vert.$$
We can now apply the Borel-Caratheodory estimate that can be found in Titchmarsh \cite{Tit}.
\begin{lem}
 Assume that $f$ is holomorphic on a neighborhood of the closed disc $\overline {D}(0,R)$, and $f(0)=0$. Then
 for all $r<R$, we have
 $$ \max_{\vert z\vert\leq r} \vert f'(z) \vert \leq \frac{8\rho_1}{(R-r)^2}\max_{\vert \zeta \vert \leq R} \vert \Re(f(z))\vert.$$
 \end{lem}
 \noindent
 Applying this lemma to $f(\zeta)=\log {\mathcal Z}_q(\zeta)$, we end up with the estimate for all 
 $\vert \zeta \vert\leq r<r'<\rho_1$,
 $$\left \vert  \frac{{\mathcal Z}'_q(\zeta)}{{\mathcal Z}_q(\zeta)} \right \vert\leq\frac{8r'}{(r'-r)^2}( 
 C_1\vert q \vert \log(1+\vert \zeta \vert)+C_2 \vert \zeta \vert)=O(\vert q \vert),$$
 the implied constant depending only on $r, \rho_1$. $\square$

\end{subsection}
\begin{subsection}{Proof in the easiest case}
We can now prove Theorem \ref{key} in the easiest case of the one dimensional torus. A first important observation
is that for all $\vert \zeta \vert$ small enough (actually $\vert \zeta \vert <1$), using the formula for the trace, the determinant ${\mathcal Z}_q(\zeta)$ is given by 
\begin{equation}
\label{trace1}
{\mathcal Z}_q(\zeta)=\exp\left(-\sum_{n=1}^{+\infty}\frac{\zeta^n}{n}\textrm{Tr}({\mathcal M}_q^n) \right)
\end{equation}
$$=\exp\left(-\sum_{n=1}^{+\infty}\frac{\zeta^n}{n} 
\sum_{T^nx=x} e^{iq\tau^{(n)}(x)} \frac{e^{\varphi^{(n)}(x)}}{1-((T^n)'(x))^{-1}}\right).$$
Given $0<r<1$, uniform convergence of the above series on the circle $\{ \vert z\vert=r \}$ yields the integration formula
\begin{equation}
\label{integral1}
-\frac{1}{2i\pi} \int_{\vert \zeta \vert=r} \frac{{\mathcal Z}'_q(\zeta)}{{\mathcal Z}_q(\zeta)} \zeta^{-n} d\zeta=
\sum_{T^nx=x} e^{iq\tau^{(n)}(x)} \frac{e^{\varphi^{(n)}(x)}}{1-((T^n)'(x))^{-1}}.
\end{equation}
For simplicity, we set for all $n\geq 1$ and $q\in \Z$, 
$$W(n,q):=\sum_{T^nx=x} e^{iq\tau^{(n)}(x)} \frac{e^{\varphi^{(n)}(x)}}{1-((T^n)'(x))^{-1}}.$$
We will use the following {\it theta inversion} \footnote{It follows from Poisson summation formula and the identity for the Fourier transform of gaussians.}
 identity valid for all $t>0$ and all $\theta \in \R$,
$$h(t,\theta):=\sum_{q\in \Z}e^{-q^2t}e^{iq\theta}=\sqrt{\frac{\pi}{t}}\sum_{p\in \Z}e^{-\frac{(\theta-2\pi p)^2}{4t}}. $$
We point out that we have for all $0<c<\pi^2/4$, as $t\rightarrow 0$, 
$$h(t,0)=\sqrt{\frac{\pi}{t}}\sum_{p\in \Z} e^{-\frac{p^2\pi^2}{t}}=\sqrt{\frac{\pi}{t}}+O\left(  e^{-c/t}\right).$$
For all $t>0$ and $n\geq 1$, we set 
$$S(t,n):=\sum_{q\in \Z}e^{-tq^2} \vert W(n,q) \vert^2,$$
and we also define for all $x\in [0,1]$ such that $T^nx=x$,
$$G_n(x):=\frac{e^{\varphi^{(n)}(x)}}{1-((T^n)'(x))^{-1} }.$$
The theta inversion formula yields 
$$S(t,n)=\sum_{T^nx=x,\ T^ny=y} G_n(x)G_n(y) \sum_{q\in \Z}e^{-q^2 t}e^{iq\tau^{(n)}(x)-iq\tau^{(n)}(y)}$$
$$=\sum_{T^nx=x,\ T^ny=y} G_n(x)G_n(y) h(t, \tau^{(n)}(x)-\tau^{(n)}(y)),$$
each term in this sum being positive.
Dropping all terms except the diagonal ones gives for all $t$ small enough, and $n$ large enough, 
$$S(t,n)\geq \frac{C_1}{\sqrt{t}}\sum_{T^nx=x} 
\frac{e^{2\varphi^{(n)}(x)}}{(1-((T^n)'(x))^{-1})^2}$$
$$\geq \frac{C_2}{\sqrt{t}}\sum_{T^nx=x} 
e^{2\varphi^{(n)}(x)},$$
for some well chosen constants $C_1,C_2>0$. Keeping this lower bound for $S(t,n)$ in mind, the end of the proof is by
contradiction. Assume that for all $\vert q \vert \geq q_0$, $q_0$ being a large integer, the spectrum of ${\mathcal M}_q$ is
included in an open disc $D(0,\rho)$ of radius $1>\rho>0$. This shows that the meromorphic function 
($\vert q\vert \geq q_0$)
$$\frac{{\mathcal Z}'_q(\zeta)}{{\mathcal Z}_q(\zeta)}$$
is actually analytic on a neighborhood of the closed disc $\overline{D}(0,\rho^{-1})$. 
Fix $\epsilon>0$.
Using Proposition \ref{bound1}, we deduce using the integral formula (\ref{integral1}) and contour deformation up to 
$(\rho+\epsilon)^{-1}$, that there exists a constant $C_\epsilon>0$ such that for all $\vert q\vert \geq q_0$,
$$\vert W(n,q)\vert \leq C_\epsilon \vert q\vert (\rho+\epsilon)^{n-1}.$$
Therefore we have the upper bound
$$S(t,n)\leq \sum_{\vert q\vert<q_0} \vert W(n,q) \vert^2 +C_\epsilon^2 (\rho+\epsilon)^{2n-2}
\sum_{\vert q\vert\geq q_0} \vert q \vert ^2 e^{-q^2t}.$$
The integral formula (\ref{integral1}) is valid for all radii $r<1$ so using proposition \ref{bound1} again we get that
$$\sum_{\vert q\vert<q_0} \vert W(n,q) \vert^2=O( (1+\epsilon)^{2n-2}).$$
On the other hand, we have using summation by parts, 
$$\sum_{ q\geq 0}  q^2 e^{-q^2t}=2\int_{1}^{+\infty}[u]u(u^2t-1)e^{-u^2t}du.$$
A change of variable $v=u\sqrt{t}$ in the integral yields as $t\rightarrow 0$, 
$$\sum_{q\geq 0} q^2 e^{-q^2 t}=O\left(  \frac{1}{t^{3/2}} \right).$$
In a nutshell, we have for $n$ large and $t$ small,
\begin{equation}
\label{end1}
C_2 \sum_{T^nx=x} 
e^{2\varphi^{(n)}(x)}\leq O( \sqrt{t} (1+\epsilon)^{2n})+O(t^{-1} (\rho+\epsilon)^{2n}),
\end{equation}
where the implied constants do not depend on $t,n$.
We recall that a classical formula on the topological pressure states that
$$\lim_{n\rightarrow +\infty} \left ( \sum_{T^nx=x} 
e^{2\varphi^{(n)}(x)}\right)^{1/n}=e^{P(2\varphi)}.$$
Now set $t=\mu^{2n}$ for some $0<\mu<1$. Going back to (\ref{end1}) and using the above formula on the pressure we get
letting $n\rightarrow +\infty$,
$$e^{P(2\varphi)}\leq \max \left \{ \mu(1+\epsilon)^2 , \left (\frac{\rho+\epsilon}{\mu} \right)^2 \right \},$$
a contradiction whenever 
$$\mu< e^{P(2\varphi)}(1+\varepsilon)^{-2}\  \textrm{and}\  \rho+\epsilon< \mu e^{\half P(2\varphi)}.$$
Since $\epsilon$ is arbitrarily small, we can clearly get a contradiction for all $\rho < e^{\frac{3}{2} P(2\varphi)}$, and the proof
is done with $\gamma(\G)=3/2$. $\square$

\bigskip
\noindent The wise reader will have noticed that $h(t,\theta)$ is the {\it heat kernel } on the torus. This will play an important role
in the proof of the general case.
\end{subsection}
\end{section}
\begin{section}{Laplacian and a priori bounds on representations}
In this section we prove the necessary upper bounds that are required to mimic the elementary proof of the previous section. 
A basic knowledge of Lie Group theory is required, our main references are \cite{Bump,Faraut}.
In the following of this paper, we assume that $\G$ is a compact connected Lie Group. We wil denote by $e_\G$ the neutral element of $\G$.
The (real) dimension of $\G$ is denoted by $d$.
A Lie Group has always a real analytic structure for which $\G$ is a real analytic manifold. Left and right translations,
inversion are real analytic diffeomorphism on $\G$. Any linear representation $\pi:\G \rightarrow \GL (V)$ is automatically real analytic. Let $\LG$ be the Lie algebra of $\G$.
Since $\G$ is compact, there exists an euclidean inner product $\beta$ on $\LG$ wich is invariant by the adjoint representation $\mathrm{ Ad}$ i.e. for all $X,Y \in \LG$, for
all $g \in \G$, 
$$\beta ( \mathrm{Ad}(g) X, \mathrm{Ad}(g) Y )=\beta(X,Y).$$
This euclidean inner product on $\LG \simeq T_{e_\G} \G$ then induces a natural left-right invariant Riemannian metric (see \cite{Bump} p.102) $\beta$ on $\G$ which is also real analytic. Clearly the Riemannian volume inherited from the metric $\beta$ is also invariant by (left-right) translations and is proportional to Haar measure $m$ so we assume for simplicity that they are equal.
Let $\Delta$ denote the (positive) Laplace Beltrami operator associated to $\beta$. Let $X_1,\ldots,X_d$ be an orthonormal basis of $\LG$ for $\beta$, then the Laplacian $\Delta$ acts on $C^\infty( G)$ by 
$$\Delta f=-\sum_{i=1}^d X_i^2(f),$$
where $X_i$ are understood as left-invariant vector fields on $\G$. Let 
$$\pi:\G\rightarrow \GL(V_\pi)$$
 be an irreducible {\it complex} representation of $\G$. Let 
$f$ be a coefficient of the representation $\pi$, i.e. for all $x\in \G$, 
$$f(x)=\mathrm {Tr}(A\pi(x)),$$
where $A \in \End (V_\pi)$. It is a classical fact (see \cite{Faraut}, p. 163) that $f$ is an eigenfunction of $\Delta$:
$$\Delta f=\kappa_\pi f,$$
where $\kappa_\pi\geq 0$ depends only on $\pi$. All the estimates will be established with respect to the spectral parameter $\kappa_\pi$. Notice that since the
eigenspaces of the Laplace operator are finite dimensional, there are only finitely many non-equivalent irreducible representations that share the same eigenvalue
$\kappa_\pi$.

Let us now use Sogge's classical bound (see \cite{Sogge}, Proposition 2.1) on $L^\infty$ norms of eigenfunctions: 
$$\Vert f \Vert_\infty\leq C \kappa_\pi^{(d-1)/4} \Vert f \Vert_{L^2(G)},$$
where $C$ depends only on $\G$ and the choice of $\beta$. In the particular case when $f(x)=\chi_\pi(x)=\mathrm{Tr}(\pi(x))$ is the {\it character} of the representation,
we get (characters of irreducible representation have a normalized $L^2$ norm)
\begin{equation}
\label{dimestimate}
\mathrm{ dim}(\pi)=\chi_\pi(e_\G)\leq C \kappa_\pi^{(d-1)/4} .
\end{equation}
Thought this dimension estimate is clearly not optimal on the $d$-torus, it is optimal for example on $\G=\mathrm{SU}_2(\C)$.
This upper bound will be useful in the rest of the paper, however optimality of the exponent is not important for us (any polynomial bound would be enough). Remark that in the examples, one can derive such a bound from the Weyl character formula. On the complex vector space $V_\pi$, there exists a hermitian inner product such that 
the representation $\pi(x)$ is unitary. We fix this hermitian structure on each $V_\pi$ and denote by $\Vert .\Vert$ the associated norm on $V_\pi$. 
Define the following norms on $\End(V_\pi)$.
$$ \Vert A \Vert_{HS}:=(\mathrm{Tr}(A^*A))^{1/2};\ \Vert A\Vert:= \sup_{\Vert v \Vert \leq 1} \Vert A v \Vert;$$
The adjoint $A^*$ is understood with respect to the above unitarizing inner product. These two norms are related by the classical inequalities
$$\Vert A \Vert \leq \Vert A \Vert_{HS}\leq \sqrt{\mathrm{dim}(\pi)} \Vert A \Vert.$$
Let us go back to the notations of $\S 3$. We recall that $\Omega$ is a complex neighborhood of $[0,1]$. The compact extension of $T$ is defined
through a skew function $\tau:[0,1]\rightarrow \G$ which is assumed to be real analytic. Since any irreducible representation $\pi$ is real analytic,
by taking $\Omega$ small enough, $z\mapsto \pi(\tau(z))$ has a holomorphic extension to $\Omega$ which is continuous on $\overline{\Omega}$ and
takes values in $\End(V_\pi)$. A priori 
\footnote{It follows from the existence of an analytic complexification $\G\subset \G_\C$, see \cite{Bump}, chapter 27, 
that every finite dimensional representation can be holomorphically extended to the complex analytic group $\G_\C$, which shows that $\Omega$
can actually be chosen uniform with respect to $\pi$. However in the sequel, we need an effective estimate and won't use the existence of $\G_\C$.}, $\Omega$ depends on $\pi$. However we will prove the following.
\begin{propo}
\label{Fefferman} There exists a small enough neighborhood $\Omega$ of $[0,1]$ such that for all irreducible complex representation $\pi$ of $\G$,
$z\mapsto \pi(\tau(z))$ has a holomorphic extension to $\Omega$ which is continuous on $\overline{\Omega}$ and enjoys the bound
$$ \sup_{z\in \Omega} \Vert \pi(\tau(z) )\Vert \leq C (\kappa_\pi)^{3(d-1)/8} e^{C\sqrt{\kappa_\pi}},$$
for some uniform $C>0$.
\end{propo}
\noindent {\it Proof}. Fix $g \in \G$ and Let $g \in \U \subset \G$ be a small coordinate neighborhood which we identify with
an open euclidean ball 
$$B(0,\rho_0)\subset \R^d \subset \C^d.$$ 
Assume that the real analytic Riemannian metric $\beta$ has a convergent power series expansion on this ball. Denote by $\Delta$ the associated
Laplace operator. Then it follows from the analysis of
Donnelly-Fefferman \cite{Fefferman} that one can find $0<\rho_1< \rho_0$ such that for each eigenfunction $F:B(0,\rho_0)\rightarrow \C$ solution of the elliptic equation
$$\Delta F= \lambda F$$
has a holomorphic extension to the (complex ball) $B_\C(0,\rho_1)\subset \C^d$ together with the bound
$$\sup_{z\in B_\C(0,\rho_1)}\vert F(z) \vert \leq e^{C\sqrt{\lambda}} \sup_{x\in B(0,\rho_0)}\vert F(x) \vert.$$
Each coefficient of the representation being an eigenfunction of the Laplacian, it is now clear that by compactness of $[0,1]$ and real analyticity of 
$\tau$, one can choose a small enough $\Omega$ which is uniform with respect to representations. To prove the desired upper bound on $\Vert \pi(\tau(z))\Vert$,
it is enough to work locally. Choose $x_0 \in [0,1]$, set $g_0:=\tau(x_0)\in \U$, where $(\U,\varphi)$ is a local chart as above with $\varphi:\U\rightarrow B(0,\rho_0)$
and $\varphi(g_0)=0$. The map $x\rightarrow \varphi(\tau(x))$ is real analytic on  a neighborhood of $x_0$. Choose $\epsilon>0$ so small that for all 
$z\in B_\C(x_0,\epsilon)$, $\varphi(\tau(z))$ is in $B_\C(0,\rho_1)$. Let $E_{ij}$, with $1\leq i,j\leq \mathrm{dim}(\pi)$ be an orthonormal basis of $\End(V_\pi)$ 
for the inner product $(A,B)\mapsto \mathrm{Tr}(A^* B)$. Each matrix coefficient $g\mapsto \mathrm{Tr}(E_{ij}\pi(g))$ is an eigenfunction of $\Delta$ so we have
for all $z\in B_\C(x_0,\epsilon)$,
$$\vert \mathrm{Tr}(E_{ij}\pi(\tau(z)))\vert \leq e^{C\sqrt{\kappa_\pi}}\sup_{g\in \G}\vert \mathrm{Tr}(E_{ij}\pi(g))
\vert$$
$$\leq e^{C\sqrt{\kappa_\pi}}\sup_{g}\Vert \pi(g) \Vert_{HS}=e^{C\sqrt{\kappa_\pi}}\sqrt{\mathrm{dim}(\pi)}$$
by Schwarz inequality and the fact that $\pi$ is unitary. Writing
$$\Vert \pi(\tau(z))\Vert_{HS}^2=\sum_{i,j}\vert \mathrm{Tr}(E_{ij}\pi(\tau(z)))    \vert^2\leq \mathrm{dim}(\pi)^3 e^{C\kappa_\pi},$$
the proof is done using the estimate (\ref{dimestimate}). $\square$
\end{section}
\begin{section}{Dynamical Zeta functions and Trace formula}
Recall that our goal is to study the point spectrum of the transfer operator 
$${\mathcal M}_\pi(f)(x):=h_\varphi^{-1}(x)\sum_{Ty=x}e^{\varphi(y)}h_\varphi(y) f(y) \pi(\tau(y)),$$ 
where $h_\varphi$ is a real analytic density fusnished by the Ruelle-Perron-Frobenius Theorem. We need to define a natural "matrix valued" function space which in our case will be
$$B^2_\pi(\Omega):=\left \{ f:\Omega \rightarrow \End(V_\pi),\:\ \textrm{holomorphic}\ 
\textrm{and}\ \int_{\Omega}\Vert f \Vert^2_{HS} dxdy<+\infty \right \}.$$
It is a Hilbert space when endowed with the norm 
$$\Vert f\Vert_{B^2}^2:=\int_{\Omega}\Vert f \Vert^2_{HS} dxdy.$$  In the sequel, we assume that $\Omega$ is a small, fixed complex neighborhood of $[0,1]$ such that
all the previous analysis holds true. We recall the notations of $\S 3$. The contracting inverse branches are denoted by $\gamma_j:\Omega\rightarrow \Omega$
and we have for all $j=1,\ldots,k$, $\overline{\gamma_j (\Omega)}\subset \Omega$. Given $z \in \Omega$, ${\mathcal M}_\pi^n(f)(z)$ can be expressed as follows
$${\mathcal M}_\pi^n(f)(z)=h_\varphi^{-1}(z)\sum_{\vert \alpha \vert=n}h_\varphi(\gamma_\alpha z)e^{\varphi^{(n)}(\gamma_\alpha z)} f(\gamma_\alpha z) \pi(\tau^{[n]}(\gamma_\alpha z)), $$
where $\alpha \in \{1,\ldots,k \}^n$ and $\gamma_\alpha:=\gamma_{\alpha_1}\circ \ldots \circ \gamma_{\alpha_n}$.
We have the following facts.
\begin{propo} 
\label{spec2}
\begin{enumerate}
\item Each ${\mathcal M}_\pi:B^2_\pi(\Omega)\rightarrow B^2_\pi(\Omega)$ is a compact trace class operator.
\item The spectral radius of ${\mathcal M}_\pi$ is bounded by $1$.
\item Let $\vert \lambda_0(\pi) \vert \geq \vert \lambda_1(\pi) \vert \geq \ldots \geq \vert \lambda_n(\pi)\vert$ denote the
sequence of eigenvalues  of ${\mathcal M}_\pi$. There exists a constant $C>0$ and $0<\rho<1$ such that for all
$\pi$ complex irreducible  and all $n\in \N$,
$$\vert \lambda_n(\pi) \vert \leq C(\kappa_\pi )^{d-1}e^{C\sqrt{\kappa_\pi}} \rho^n.$$
\item For all $\pi$ complex irreducible and all $n\in \N$, we have the formula for the trace of ${\mathcal M}_\pi^n$
$$\textrm{Tr}({\mathcal M}_\pi^n)=\mathrm{dim}(\pi)\sum_{T^nx=x} \chi_\pi(\tau^{[n]}(x)) \frac{e^{\varphi^{(n)}(x)}}{1-((T^n)'(x))^{-1}},$$
where $\chi_\pi(g)$ is the character of $\pi$. The sum runs over $n$-periodic points of $T$ and $\tau^{[n]}(x)=\tau(x)\ldots \tau(T^{n-1}x)$.
\end{enumerate}
\end{propo}
\noindent
{\it Proof}. $(1)$ We let the reader check that ${\mathcal M}_\pi$ is indeed bounded and compact by following the same ideas as in the scalar case. The fact that
it is trace class will follow from the summability of the singular values estimates.

\noindent $(2)$ It is enough to estimate the spectral radius when acting on continuous functions $C^0(I,\End(V_\pi))$. Indeed, any $B^2_\pi(\Omega)$ eigenvalue
is also a $C^0(I,\End(V_\pi))$ eigenvalue. We endow the space $C^0(I,\End(V_\pi))$ with the obvious norm
$$ \Vert f \Vert_0:=\sup_{x\in I} \Vert f(x) \Vert.$$ 
Let $f \in C^0(I,\End(V_\pi))$ and write (using unitarity of $\pi$)
$$ \Vert {\mathcal M}_\pi^n(f)(x) \Vert\leq 
 h_\varphi^{-1}(x)\sum_{\vert \alpha \vert=n}h_\varphi(\gamma_\alpha x)e^{\varphi^{(n)}(\gamma_\alpha x)} \Vert f(\gamma_\alpha x)\Vert,$$
 $$\leq C \left ( \sum_{\vert \alpha \vert=n} e^{\sup_{x \in I} \varphi^{(n)}(\gamma_\alpha x)} \right ) \Vert f \Vert_0.$$
 Now the claim follows readily from the formula of the topological pressure:
 $$\lim_{n\rightarrow +\infty} \left ( \sum_{T^nx=x} 
e^{2\varphi^{(n)}(x)}\right)^{1/n}=e^{P(\varphi)}=1,$$
and the fact that there exists $M>0$ such that for all $n$ and all $x\in I$, all $\alpha \in \{1,\ldots,k\}^n$,
$$\left \vert \frac{d}{dx} \varphi^{(n)}(\gamma_\alpha x) \right \vert \leq M$$ which is a standard fact in uniformly hyperbolic dynamics
referred as "bounded distortion property".

\noindent $(3)$ This estimate will follow from a singular value estimate. We recall that the $n$-th singular value of the compact operator ${\mathcal M}_\pi$ is by definition the $n$-th eigenvalue (ordered decreasingly) of the positive compact self-adjoint operator  $\sqrt{ {\mathcal M}_\pi^* {\mathcal M}_\pi}$. They are given by Courant's minimax
formula 
$$s_n({\mathcal M}_\pi)=\min_{\mathrm{dim}(F)=n} \max_{f\in F^\perp, \Vert f \Vert \leq 1} \Vert {\mathcal M}_\pi(f) \Vert_{\B^2}.$$ 
In particular, we can derive the fact that for any basis $(e_l)_{l \in \N}$ of $B^2_\pi(\Omega)$, 
$$s_n({\mathcal M}_\pi)\leq \sum_{l\geq n} \Vert {\mathcal M}_\pi(e_l)\Vert_{B^2}.$$
To be able to use that estimate, we need an "explicit" basis of $B^2(\Omega)$. We point out that if $\Omega$ is of type $\Omega_\epsilon:=I+B(0,\epsilon)$ then
the boundary is actually "Dini-smooth" (see \cite{Pommerenke}) and there exists a conformal mapping $\psi:\Omega \rightarrow {\mathbb D}$, where $\mathbb D$ is the unit
disc such that both $\psi$ and $\psi'$ are continuous up to the boundary of $\Omega$. There is thus an isometry 
$${\mathcal I}: \left \{ \begin{array}{ccc}
B^2({\mathbb D})& \rightarrow & B^2(\D)\\
f& \mapsto & \psi' (f\circ \psi).
\end{array} \right.$$
Viewing our Hilbert space $B^2_\pi(\Omega)$ as
$$B^2_\pi(\Omega)=\bigoplus_{1\leq i,j\leq \mathrm{dim}(\pi)} B^2(\Omega),$$
we have a natural Hilbert basis given by the family ($l\in \N$, $1\leq i,j\leq \mathrm{dim}(\pi)$)
$${\mathbf e}_l^{ij}(z):=\left(\frac{l+1}{\pi}\right)E_{ij}\psi'(z) (\psi(z))^l,$$
Where $E_{ij}$ is the linear map of $V_\pi$ into itself whose matrix in a fixed orthonormal basis of $V_\pi$ has zero entries everywhere except
$1$ for row $i$ and column $j$.

For all $z\in \Omega$, we write
$$\Vert {\mathcal M}_\pi( \mathbf{e}^{ij}_l)(z)\Vert_{HS}$$
$$\leq C_1 k \sup_{y\in \Omega}e^{\varphi(y)}\sqrt{\mathrm{dim}(\pi)} (\kappa_\pi)^{3\frac{d-1}{8}} e^{C\sqrt{\kappa_\pi}}
\max_{m}\Vert  \mathbf{e}^{ij}_l(\gamma_m z) \Vert_{HS}, $$
But we have 
$$\Vert \mathbf{e}^{ij}_l(\gamma_m z) \Vert_{HS}^2=\left \vert  \left(\frac{l+1}{\pi}\right)\psi'(\gamma_m z) (\psi(\gamma_m z))^l\right \vert^2 $$
and therefore
$$\Vert \mathbf{e}^{ij}_l(\gamma_m z) \Vert_{HS}^2\leq C_2 l^2 \rho^{2l},$$ 
where 
$$\rho=\sup_{z \in \Omega,\ 1\leq m\leq k} \vert \psi(\gamma_m z)\vert <1.$$
To conclude this estimate, we have
$$ s_n({\mathcal M}_\pi)\leq \sum_{l\geq n} \sum_{i,j} \Vert \mathbf{e}^{ij}_l(\gamma_m z) \Vert_{B^2}$$
$$\leq C_3 \kappa_{\pi}^{d-1} e^{C\sqrt{\kappa_\pi}} (\widetilde{\rho})^n,$$
where $\rho<\widetilde{\rho}<1$. To relate it with the eigenvalue estimate, we simply write
$$ \vert \lambda_n\vert \leq \left (\prod_{k=1}^n \vert \lambda_k \vert \right)^{1/n}\leq\left (\prod_{k=1}^n \vert s_k \vert \right)^{1/n}, $$
where we have used (in the second equality) a standard Weyl inequality (see \cite{Simon} on Weyl inequalities) which gives the desired estimate.

\noindent $(4)$ The trace is usually defined by the following formula:
$$ \mathrm{Tr}({\mathcal M}_\pi^n):=\sum_{l=0}^{+\infty} \sum_{i,j} \langle \mathbf{e}^{ij}_l, {\mathcal M}_\pi^n(\mathbf{e}^{ij}_l  )\rangle_{B^2}.$$
Let us write
$$\sum_{i,j} \langle \mathbf{e}^{ij}_l, {\mathcal M}_\pi^n(\mathbf{e}^{ij}_l  )\rangle_{B^2}=$$
$$\sum_{\vert \alpha \vert=n} \sum_{i,j} \int_{\Omega} h^{-1}_\varphi(z) h_\varphi(\gamma_\alpha z) e^{\varphi^{(n)}(\gamma_\alpha z)}
\overline{{\mathbf e}_l(z)}{\mathbf e}_l(\gamma_\alpha z) \Psi_{i,j,\alpha}(z) dxdy,$$
where we have set
$$\Psi_{i,j,\alpha}(z):=\mathrm{Tr} \left(E_{ij}^* E_{ij} \pi(\tau^{[m]}(\gamma_\alpha z))  \right). $$
Since we have
$$E_{ij}^* E_{ij}=E_{jj},$$
we have $$ \sum_{i,j}\Psi_{i,j,\alpha}(z)=\mathrm{dim}(\pi) \mathrm{Tr}( \pi(\tau^{[m]}(\gamma_\alpha z))).    $$
We end up with 
$$\sum_{i,j} \langle \mathbf{e}^{ij}_l, {\mathcal M}_\pi^n(\mathbf{e}^{ij}_l  )\rangle_{B^2}=$$
$$\sum_{\vert \alpha \vert=n} \mathrm{dim}(\pi) \int_{\Omega} h^{-1}_\varphi(z) h_\varphi(\gamma_\alpha z) e^{\varphi^{(n)}(\gamma_\alpha z)}
\overline{{\mathbf e}_l(z)}{\mathbf e}_l(\gamma_\alpha z) \chi_\pi(\tau^{[n]}(\gamma_\alpha z)) dxdy.$$
Summing over all $l$, using uniform convergence and Fubini, we are left to compute
$$\mathrm{Tr}({\mathcal M}_\pi^n)=\sum_{\vert \alpha \vert=n} \int_{\Omega} F_n(z) B_\Omega(\gamma_\alpha z, z) dxdy,$$
where $ F_n(z)=\mathrm{dim}(\pi) h^{-1}_\varphi(z) h_\varphi(\gamma_\alpha z) e^{\varphi^{(n)}(\gamma_\alpha z)}\chi_\pi(\tau^{[n]}(\gamma_\alpha z))$,
and $B_\Omega(z,\zeta)$ is the classical Bergman kernel of the domain $\Omega$. Using the conformal change of variable,
we can assume that $\Omega$ is the unit disc and it is now a standard computation (using Stokes theorem) that yields the desired formula. $\square$

We can now state the analog of Proposition  \ref{bound1}. We set in the following for all $\zeta \in \C$,
$${\mathcal Z}_\pi(\zeta):=\det(I-\zeta \mathcal{M}_\pi).$$
\begin{propo}
\label{bound2}
Assume that for all $\kappa_\pi$ large, ${\mathcal Z}_\pi(\zeta)$ does not vanish on the open disc $D(0,\rho_1)$, for some $\rho_1>0$. Then for all $r<\rho_1$, one can find $C_r>0$ such that for all $\kappa_\pi$ large and all $\zeta \in D(0,r)$,
$$\left \vert  \frac{{\mathcal Z}'_\pi(\zeta)}{{\mathcal Z}_\pi(\zeta)} \right \vert\leq  C_r \sqrt{\kappa_\pi}.$$
\end{propo}
\noindent {\it Proof}. Use the bound $(2)$ and $(3)$ from the previous proposition and mimic the proof of Proposition \ref{bound1}. One can observe
that the polynomial terms in estimate $(3)$ do not contribute to the final estimate. $\square$
\end{section}
\begin{section}{Heat kernels on Compact Lie groups}
Our basic reference for that section is \cite{Davies}.
Let $M$ be a compact $d$-dimensional Riemannian Manifold. let $\Delta$ be the (positive) Laplacian on $M$. Let $f\in L^2(M)$ and consider the
Cauchy problem (heat equation):
\begin{equation}
\label{heat} 
\left \{ \begin{array}{ccc}
\Delta u(x,t)+\partial_t u(x,t)&=&0,\ t>0\\
  \lim_{t\rightarrow 0} u(x,t)&=&f(x) \in L^2.
\end{array}    \right.
\end{equation}
It is a well known fact, see \cite{Davies}, chapter 5, that there exists a unique kernel $H(t;x,y)\in C^\infty( (0,+\infty)\times M\times M)$ such that
for $t>0$,
$$u(t,x)=\int_M H(t;x,y)f(y) d\mathrm{Vol}(y),$$
solves the problem (\ref{heat}). Moreover $H(t;x,y)>0$ for all $t>0$. A basic lower bound which will be useful for us is following. One can find a constant
$C>0$ such that for all $t>0$ small enough and $x,y\in M$,
$$H(t;x,y)\geq \frac{C}{t^{d/2}} e^{-\frac{\mathrm{dist}(x,y)^2}{4t}},$$
where $\mathrm{dist}(x,y)$ is the Riemannian distance. This estimate follows from a Harnack estimate of Li and Yau (see \cite{Davies}, P. 162)
and the local asymptotics on the diagonal as $t\rightarrow 0$. This estimate is actually valid on non compact complete manifolds with Ricci curvature
bounded from below. If we assume that $M=\G$ has a Lie group structure then a formula for the heat kernel $H(t;x,y)$ can be derived using
characters of representations. More precisely, we have the following.
\begin{propo}
\label{heatformula}
Assume that $\G$ is a compact Lie Group endowed with a left-right invariant Riemannian metric. Then for all $t>0$, we have
$$H(t;x,y)=h_t(xy^{-1}):=\sum_{\pi} e^{-t\kappa_\pi}\mathrm{dim}(\pi)\chi_\pi(xy^{-1}),$$
the sum (on all irreducible representations) being uniformly convergent.
\end{propo}
\noindent {\it Proof}. This result is fairly standard (see \cite{Faraut}, chapter 12), however we will encounter in the proof some estimates that will be
important later.
 First we recall that by Weyl asymptotics, as $R\rightarrow +\infty$,
$$ \# \{ \lambda \in \mathrm{Spec}(\Delta)\ :\ \lambda \leq R\}\sim C_d R^{d/2},$$
eigenvales being counted with multiplicities. Since the space ${\mathit M}_\pi$ of coefficients of the representation $\pi$ has dimension $\mathrm{dim}(\pi)^2$ and each $F\in {\mathit M}_\pi$ satisfies $ \Delta F=\kappa_\pi F$, we deduce that
$$N(R):=\sum_{\kappa_\pi\leq R} \mathrm{dim}(\pi)^2 \leq C R^{d/2}.$$
Let $f \in L^2(G)$ and set for all $t\geq 0$,
$$u(x,t):=\sum_{\pi} e^{-\kappa_\pi t} \mathrm{dim}(\pi) \mathrm{Tr}( \widehat{f}(\pi) \pi(x)),$$
where $$\widehat{f}(\pi)=\int_G f(g) \pi(g^{-1}) dm(g).$$
By Peter-Weyl theorem (see for example \cite{Faraut}, Chapter 6), this above expression makes sense in $L^2(G)$ and is equal to $f$ for $t=0$.
Let us show that $u(x,t)$ indeed solves the heat equation. First remark that by Schwarz inequality (twice), we have
$$\vert \mathrm{Tr}( \widehat{f}(\pi)\pi(x))\vert\leq \Vert \widehat{f}(\pi) \Vert_{HS} \Vert \pi(x) \Vert_{HS} \leq \mathrm{dim}(\pi) \Vert f \Vert_{L^2(G)}.$$
This allows us to bound using Stieltjes integrals ($t>0$)
$$ \vert u(x,t) \vert \leq \sum_{\pi} e^{-\kappa_\pi t} \mathrm{dim}(\pi)^2 \Vert f \Vert_{L^2} =\Vert f \Vert_{L^2} \int_0^{+\infty} e^{-vt}d N(v).$$
Using summation by parts, the above integral is clearly seen to be convergent since $N(v)=O(v^{d/2})$. Therefore the above series are uniformly
convergent for all $t>0$.
Similarly, we have 
$$ \vert \Delta (\mathrm{Tr}( \widehat{f}(\pi)\pi(x)) ) \vert \leq \kappa_\pi \mathrm{dim}(\pi) \Vert f \Vert_{L^2(G)},$$
so that 
$$ \sum_{\pi} e^{-\kappa_\pi t} \mathrm{dim}(\pi) \Delta \left (\mathrm{Tr}( \widehat{f}(\pi) \pi(x)) \right)$$
is uniformly convergent, allowing us to write 
$$ \Delta u(x,t)= \sum_{\pi} e^{-\kappa_\pi t} \mathrm{dim}(\pi) \kappa_\pi \mathrm{Tr}( \widehat{f}(\pi) \pi(x)). $$
On the other hand, a similar argument shows that
$$\partial_t u(x,t)=-\sum_{\pi} e^{-\kappa_\pi t} \mathrm{dim}(\pi) \kappa_\pi \mathrm{Tr}( \widehat{f}(\pi) \pi(x)).$$
It remains to check that $\lim_{t\rightarrow 0}u(x,t)=f(x)$ in $L^2(G)$ sense.  This follows easily from Plancherel formula
$$\Vert u(x,t)-f(x) \Vert^2_{L^2}=\sum_\pi \mathrm{dim}(\pi)\Vert \widehat {f}(\pi)(1-e^{-\kappa_\pi t})\Vert^2_{HS}  $$
$$=(1-e^{-\kappa_\pi t})^2 \Vert f\Vert_{L^2}^2. $$
Using uniform convergence, we can finally write  $(t>0)$
$$u(x,t)=\sum_{\pi} e^{-\kappa_\pi t} \mathrm{dim}(\pi) \mathrm{Tr}\left(  \int_G f(g) \pi(g^{-1}) dm(g)\pi(x) \right)$$
$$=\sum_{\pi} e^{-\kappa_\pi t} \mathrm{dim}(\pi) \int_G f(g) \chi_\pi(g^{-1}x) dm(g)  $$
$$=\int_{G}h_t(xg^{-1})f(g)dm(g).$$
Uniqueness of the heat kernel concludes the proof. $\square$

\end{section}
\begin{section}{Proof in the abstract case} 
The proof starts with the integral formula ($r>1$)
\begin{equation}
\label{integral2}
-\frac{1}{2i\pi} \int_{\vert \zeta \vert=r} \frac{{\mathcal Z}'_\pi(\zeta)}{{\mathcal Z}_\pi(\zeta)} \zeta^{-n} d\zeta=
\mathrm{dim}(\pi)\sum_{T^nx=x} \chi_\pi(\tau^{[n]}(x))\frac{e^{\varphi^{(n)}(x)}}{1-((T^n)'(x))^{-1}}.
\end{equation}
We set for all $\pi$ irreducible and $n\geq 1$
$$W(n,\pi):=\mathrm{dim}(\pi)\sum_{T^nx=x} \chi_\pi(\tau^{[n]}(x))\frac{e^{\varphi^{(n)}(x)}}{1-((T^n)'(x))^{-1}}.$$
The first step is to examine the following average on representations: ($t>0$)
$$\mathcal{S}(n,t):=\sum_{\pi } e^{-\kappa_\pi t} \frac{\vert W(n,\pi)\vert^2}{\mathrm{dim}(\pi)^2}.$$
\begin{propo}
 There exists a constant $A>0$ such that for all $t>0$ small and $n\geq 1$,
 $$\mathcal{S}(n,t)\geq A t^{-\frac{\mathrm{rank}(\G)}{2}}\sum_{T^n x=x}e^{2\varphi^{(n)}(x)}, $$
 where $\mathrm{rank}(\G)$ is the real dimension of maximal torii.
\end{propo}
\noindent
{\it Proof.} By a direct computation, we have 
$$ \mathcal{S}(n,t)=\sum_{T^nx=x \atop T^n x'=x'}\frac{e^{\varphi^{(n)}(x)}}{1-((T^n)'(x))^{-1}}\frac{e^{\varphi^{(n)}(x')}}{1-((T^n)'(x'))^{-1}} 
F(t;\tau^{[n]}(x),\tau^{[n]}(x')),$$
where we have set
$$F(t;a,b):=\sum_{\pi} e^{-t\kappa_\pi}\chi_\pi(a)\overline{\chi_\pi(b)}.$$
Using the formula
$$\int_G\chi_\pi(xgyg^{-1})dm(g)=\frac{\chi_\pi(x)\chi_\pi(y)}{\mathrm{dim}(\pi)},  $$
and uniform convergence of the above series for $t>0$, we get
$$F(t;a,b)=\int_G h_t( agb^{-1}g^{-1}) dm(g),$$
where $h_t(x)$ is the heat kernel, see previous $\S$. The heat kernel being positive, we drop all the off-diagonal terms to get
$$\mathcal{S}(n,t)\geq \sum_{T^nx=x }\frac{e^{2\varphi^{(n)}(x)}}{(1-((T^n)'(x))^{-1})^2}
F(t;\tau^{[n]}(x),\tau^{[n]}(x)).$$
We therefore need to estimate from below the quantity $F(t;a,a)$ for all $a\in \G$. First observe that if $\G$ is commutative i.e. is a $d$-torus,
then $F(t,a,a)$ is just $h_t(e_\G)$ where $e_\G$ is the neutral element, and the claim follows readily using the short time asymptotics of the heat kernel.
If $\G$ is not a torus, then we can still use Cartan's theorem (see \cite{Bump}) which says the following things.
Call a maximal torus a maximal (for the dimension) abelian subgroup of $\G$. 
\begin{enumerate}
\item All maximal torii are conjugated.
\item Every $g\in \G$ is inside a maximal torus.
\end{enumerate}
Given a maximal torus $\mathbf{T}\subset \G$ we detone by $\mathbf{T_\epsilon}$ the $\epsilon$-tube around $\mathbf{T}$ defined by
$$\mathbf{T}_\epsilon:=\{ g \in \G\ :\ \mathrm{dist}(g,\mathbf T)\leq \epsilon \} ,$$
where "$\mathrm{dist}$" is the left-right invariant Riemannian distance. Notice that by property $(1)$ of Cartan's theorem and invariance of both distance and Haar volume, for all maximal
torii $\mathbf{T,T'}$, we have for all $\epsilon>0$
$$m(\mathbf{T}_\epsilon)=m(\mathbf{T}_\epsilon'),$$
where $m$ is Haar measure. Given $a\in \G$, pick a maximal torus $\mathbf{T}\ni a$
 and write 
 $$ F(t;a,a)\geq \int_{\mathbf T_\epsilon} h_t(aga^{-1}g^{-1})dm(g).$$
Using the Gaussian lower bound for the heat kernel from $\S 6$, we have therefore
$$ F(t;a,a)\geq \frac{C}{t^{d/2}}\int_{\mathbf T_\epsilon} \exp \left ( -\frac{\mathrm{dist}^2(aga^{-1}g^{-1},e_\G)}{4t} \right)dm(g).$$
Given $g \in \mathbf{T}_\epsilon$, pick $\widetilde{g} \in \mathbf{T}$ such that $\mathrm{dist}(g,\widetilde{g})\leq 2\epsilon$. We have
$$\mathrm{dist}(aga^{-1}g^{-1},e_\G)=\mathrm{dist}(ag,ga)$$
$$\leq \mathrm{dist}(ag,a\widetilde{g})+\underbrace{\mathrm{dist}(a\widetilde{g},\widetilde{g}a)}_{=0}+\mathrm{dist}(\widetilde{g}a,ga)$$ 
$$=2\mathrm{dist}(g,\widetilde{g})\leq 4\epsilon.$$
Taking $\epsilon=\sqrt{t}$, we get for all $t>0$ small,
\begin{equation}
\label{heatvolume}
F(t;a,a)\geq \frac{Ce^{-4}}{t^{d/2}} m(\mathbf{T}_{\epsilon}).
\end{equation}
Remark now that $m(\mathbf{T}_{\epsilon})$ does not depend on the chosen maximal torus (and hence $a$) and we can fix a maximal torus $\mathbf{T}$ once for all to compute 
$m(\mathbf{T}_{\epsilon})$. Since Haar mesure is proportional to Riemannian volume and maximal torii are smooth submanifolds, 
the final estimate will follow from an easy volume estimate (from below). Set for simplicity $r=\textrm{rank}(\G)$. 
By $B_n(0,\alpha)$ we mean the euclidean ball of $\R^n$ of radius $\alpha$ centered at $0$. The usual euclidean norm on $\R^n$
is denoted by $\vert x \vert$.
Let $\U$ be a coordinate
neighborhood around $e_\G$ and $\psi:\U\rightarrow \R^d$ a coordinate map such that 
$$\psi(\U)=B_r(0,\alpha_1)\times B_{d-r}(0,\alpha_2).$$
Since $\mathbf T$ is a submanifold, we can linearize it locally and choose $\psi$ so that 
$$\psi(\mathbf{ T}\cap \U)=B_r(0,\alpha_1)\times \{ 0 \}.$$
Without loss of generality, we can assume that the Riemannian metric $g=\sum g_{ij}dx^idx^j$ in this coordinate patch is continuous up to
the boundary. Let $d_g(x,y)$ denote the associated distance on $\psi(\U)$. Compactness then shows that for all $x,y \in \psi(\U)$,
$$d_g(x,y)\leq C_1 \vert x-y \vert,$$
for some $C_1>0$. For all $\epsilon>0$ small enough, we have 
$$\psi(\mathbf{T}_\epsilon\cap \U)\supset B_r(0,\alpha_1)\times B_{d-r}(0, \epsilon/C_1).$$
Therefore, 
$$m(\mathbf{ T}_\epsilon)\geq m(\mathbf{T}_\epsilon\cap \U)=\int_{\psi(\mathbf{T}_\epsilon\cap \U)} \sqrt{\det( g_{ij})}dx_1\ldots dx_d$$
$$\geq C_2 \mathrm{Leb}(B_r(0,\alpha_1)\times B_{d-r}(0, \epsilon/C_1))=C_3 \epsilon^{d-r},$$
here $\mathrm{Leb}$ is of course Lebesgue measure. Going back to formula (\ref{heatvolume}) and writing $\epsilon=\sqrt{t}$,
the proof is done. $\square$

\bigskip
\noindent
We can now go back to the main proof of Theorem \ref{key}. Fix $R>0$ large and assume that for all irreducible $\pi$ with $\kappa_\pi \geq R$,
the compact operator $\mathcal {M}_\pi: B^2_\pi(\Omega)\rightarrow B^2_\pi(\Omega)$ has no eigenvalues outside the disc
$$\{ z \in \C\ :\ \vert z \vert < \rho<1 \}.$$ 
Each meromorphic function $\zeta\mapsto \mathcal{Z}'_\pi(\zeta)/\mathcal{Z}_\pi(\zeta)$ is therefore analytic inside 
$$\{ \zeta\in \C\ :\ \vert \zeta\vert \leq \rho^{-1}\}.$$
Using the integral formula (\ref{integral2}), we can perform a contour shift and use the upper bound of Proposition (\ref{bound2}) to deduce that for $\kappa_\pi\geq R$,
we have for all $\epsilon>0$,
$$\vert W(n,\pi)\vert \leq C_\epsilon(\rho+\epsilon)^{n-1}\sqrt{\kappa_\pi}. $$
On the other hand, for all $\pi$ with $\kappa_\pi<R$ (which corresponds to a finite number of non-equivalent irreducible representations),
we can use the trace formula and the fact that $P(\varphi)=0$ to write for all $\epsilon>0$,
$$\vert W(n,q)\vert \leq C_\epsilon (1+\epsilon)^n.$$
Gathering our estimates, we have for all $\epsilon>0$, all $t>0$ small and $n\geq 1$,
$$A t^{-\mathrm{rank}(\G)/2} \sum_{T^n x=x}e^{2\varphi^{(n)}(x)}\leq O((1+\epsilon)^{2n})$$
$$+O\left(  (\rho+\epsilon)^{2n}
\sum_{\kappa_\pi\geq R} e^{-t\kappa_\pi}\frac{\kappa_\pi}{(\mathrm{dim}(\pi))^2}    \right).$$
Assume we have established that as $t\rightarrow 0$, 
$$\sum_{\pi} e^{-t\kappa_\pi}\frac{\kappa_\pi}{(\mathrm{dim}(\pi))^2}=O\left(  t^{-\beta} \right),$$
for some $\beta>0$. Choosing $t=e^{-\alpha n}$ we get
$$e^{\alpha \mathrm{rank}(\G)/2}\left ( \sum_{T^n x=x}e^{2\varphi^{(n)}(x)}\right)^{1/n}$$
$$\leq \left ( O((1+\epsilon)^{2n}) +O(  e^{\alpha\beta n}(\rho+\epsilon)^{2n}) \right)^{1/n}.$$
Taking the limit as $n \rightarrow +\infty$, we have
$$e^{\alpha \mathrm{rank}(\G)/2} e^{P(2\varphi)}\leq \max \{(1+\epsilon)^2 ,(\rho+\epsilon)^2e^{\alpha \beta}\}.$$
We have a contradiction if we take
$$\alpha=\frac{2\vert P(2\varphi)\vert}{\mathrm{rank}(\G)}+3\epsilon,$$
and 
$$ \rho+\epsilon<(1+\epsilon) e^{-\alpha\beta/2}=(1+\epsilon)e^{-3\epsilon \beta/2}e^{\beta \frac{P(2\varphi)}{\mathrm{rank}(\G)}}.$$
We can therefore take for $\gamma(\G)=\beta/\mathrm{rank}(\G)$. 
It remains to compute $\beta$ in general.  The quality of the estimate highly depends on the informations that we have
on $\mathrm{dim}(\pi)$.  Clearly if we have no available {\it lower bound} for $\mathrm{dim}(\pi)$, we can write 
$$\sum_{\pi} e^{-t\kappa_\pi}\frac{\kappa_\pi}{(\mathrm{dim}(\pi))^2}\leq \sum_{\pi} e^{-t\kappa_\pi}\kappa_\pi $$
$$\leq \int_0^{+\infty }e^{-t u} u dN(u),$$
where $N(u)$ is the counting function for the eigenvalues defined in $\S 6$. 
Summation by parts yields
$$\int_0^{+\infty }e^{-t u} u dN(u)=O(1)+\int_0^{+\infty}e^{-tu}(tu-1)N(u)du$$
$$=O\left( t\int_0^{+\infty} e^{-tu} u^{d/2+1}du \right)=O\left( t^{-(1+d/2)}\right). $$
This leads to the value $\gamma(\G)=\frac{1+d/2}{\mathrm{rank}(\G)}.$
Below we list the values of $\gamma(\G)$ for the most classical groups.

\begin{center}
\begin{tabular}{|c|c|}
  \hline
  Group $\G$ & Constant $\gamma(\G)$ \\
  \hline
  
   $\mathbf{T}^d=\R^d/ 2\pi \Z^d$ & $\scriptstyle{\frac{1}{2}+\frac{1}{d}}$    \\
   
       $U_n(\C)$ &$\scriptstyle{\frac{n}{2}+\frac{1}{n}} $       \\
       $SU_{2p}(\C)$& $\scriptstyle{\frac{1}{2p}+2p}$\\
       $SU_{2p+1}(\C)$&  $\scriptstyle{\frac{1}{p}+2p+2}$\\
       $SO_{2p}(\R)$&  $\scriptstyle{\frac{1}{p}+\frac{2p-1}{2}}$\\
       $SO_{2p+1}(\R)$& $\scriptstyle{\frac{1}{p}+\frac{2p+1}{2}}$\\
  \hline
\end{tabular}
\end{center}
\noindent
We point out that except in the case of torii where $\mathrm{dim}(\pi)=1$, the values of $\gamma(\G)$ can be improved when some explicit enough formulas are valid, which should be the case of some of the
groups listed above. Let us just have a look at the example of $SU_2(\C)$.  In that case the irreducible complex representations $\pi_m$ are
parametrized by $m\in \N$ and we have $\kappa_{\pi_m}=m(m+2)$, $\mathrm{dim}(\pi_m)=m+1$, see for example \cite{Faraut}. 
Going back to the above estimate, we are left with
$$\sum_{\pi} e^{-t\kappa_\pi}\frac{\kappa_\pi}{(\mathrm{dim}(\pi))^2}=O\left( \sum_{m\in \N} e^{-tm(m+2)} \right)$$
$$=O\left( t^{-1/2}\right).$$
We have obtained $\gamma(SU_2(\C))=1/2$ instead of the $5/2$ given by the above formula. Being
doubly covered by $SU_2(\C)$, the same estimate holds for $SO_3(\R)$.

\end{section}

\end{document}